\input amstex
\documentstyle{amsppt}
\magnification=\magstep1
\NoRunningHeads

\vcorrection{-1cm}

\topmatter
\title  A SHORT PROOF OF SOLUTION FORMULAS FOR \\
THE LINEAR DIFFERENTIAL EQUATIONS \\ WITH CONSTANT COEFFICIENTS
\endtitle
\author Evgeniy Pustylnik \\\\ Technion, Haifa 32000, Israel \\
{\it E-mail: evg\@ technion.ac.il}
\endauthor
\abstract
The solution of equations from the title is well known since the Euler's
time. However, its proof in the case of multiple roots of the characteristic
polynomial is rather long and technical and even appearance of the factors
$x^m$ looks artificial. A new proof, proposed below, seems not only shorter,
but also more comprehensible for the students of any level.
\endabstract
\endtopmatter

\document

{\bf 1.} Most textbooks, devoted to ordinary differential equations, have 
special sections on the linear equations with constant coefficients (see, 
e.g., very popular books [1], [2] etc.). The methods of proofs are classical 
and similar everywhere, mainly suited to engineering (non-mathematician) 
students. However, in the case of multiple roots of the corresponding 
characteristic  polynomials, the proofs become either too long (if detailed 
enough) or too concise, requiring great own efforts of readers. Consequently, 
many lecturers simply omit the proofs, making all formulas rather mysterious 
for the students. In the present note we try to avoid this obstacle. \medskip

{\bf 2.} Let $L=L[y]$ be a linear differential operator with constant
coefficients
$$ L[y]=y^{(n)}+p_1y^{(n-1)}+\ldots+p_{n-1}y'+p_ny $$
and let $P=r^n+p_1r^{n-1}+\ldots+p_{n-1}r+p_n$ mean its characteristic
polynomial. When the operators have the subscripts $L_1, L_2,\ldots$, the
same subscripts will be used for the corresponding polynomials.

Obviously, all operators as well as their characteristic polynomials
form isomorphic linear spaces. Moreover, these spaces are isomorphic
commutative algebras if the product of operators is understood as
their composition: $L_1L_2=L_1(L_2)=L_2(L_1)$ (the commutativity is
due to the constant coefficients). As a result, the standard
decomposition of a characteristic polynomial
$$ P(r)=(r-r_1)^{m_1}(r-r_2)^{m_2}\cdots(r-r_k)^{m_k}, \tag1 $$
where all the numbers $r_1,r_2\,,\ldots,r_k$ are different, implies
an analogous decomposition $L=L_1L_2\ldots L_k$, in which every $L_i\,,\
i=1\,,\ldots,k$, has the simplest characteristic polynomial $P_i(r)=
(r-r_i)^{m_i}$.

\proclaim{Lemma 1} Let an operator $L$ have the above mentioned decomposition
and let a function $y(x)$ be such that $L_i[y]=0$ for some $i=1\,,\ldots,k$.
Then $y(x)$ is a solution of the whole equation $L[y]=0$. \endproclaim

\demo{Proof} The commutativity of operators allows us to move $L_i[y]$ to any
place in the product, e.g., to assume that $i=k$. Then
$$ L[y]=L_1L_2\ldots L_{k-1}L_k[y]=L_1L_2\ldots L_{k-1}[0]=0. $$ \enddemo

Now we proceed to solution of equations with the simplest characteristic
polynomials. For convenience, we temporarily omit the subscripts in notation
of operators.

\proclaim{Lemma 2} Let $L$ be a differential operator with the characteristic
polynomial $P(r)=(r-a)^m$. Then the general solution of the equation $L[y]=0$ is
$$ y(x)=(C_1+C_2x+\ldots+C_mx^{m-1}\!)\,e^{ax}, \tag2 $$
where $C_1, C_2, \ldots, C_m$ are arbitrary constants.  \endproclaim

\demo{Proof} Let us write $P(r)$ with binomial coefficients
$$ P(r)=\sum_{i=0}^m \binom mi r^{m-i}(-a)^i, $$ then the differential equation will be
$$ L[y]=\sum_{i=0}^m \binom mi (-a)^iy^{(m-i)}=0. $$ Now we multiply both sides of this
equation by the function $v=e^{-ax}$ so that $v^{(i)}=(-a)^ie^{-ax}$. Thus the differential
equation obtains the form $$ \sum_{i=0}^m \binom mi v^{(i)}y^{(m-i)}=0. $$
But the left-hand side of this equality is exactly the Leibnitz expression for the derivative
$(vy)^{(m)}$, so that $(vy)^{(m)}=0$. This implies that $vy$ is an arbitrary polynomial of
the degree $m-1$, namely, $vy=C_1+C_2x+\ldots+C_mx^{m-1}$. Dividing by $v=e^{-ax}$, we get
the proof of Lemma. \enddemo

If we take now $a=r_i,\ m=m_i,\ i=1\,,\ldots,k$, given in the decomposition (1), we obtain from
(2) the general solutions $y_i(x)$ for all equations $L_i[y]=0$, discussed in Lemma 1. By the
same Lemma, all $y_i(x)$ are solutions of the equation $L[y]=0$ as well as any linear combination
of them. Thus we arrive at the following main assertion.

\proclaim{Theorem} Let a linear differential operator $L[y]$ with constant
coefficients has the characteristic polynomial $P(r)$ with the decomposition (1).
Then the general solution of the differential equation $L[y]=0$ has a form
$$ y(x)=\sum_{i=1}^k(C_{i1}+C_{i2}x+\ldots C_{im_i}x^{m_i-1})e^{r_ix}, \tag3 $$
where all $C_{ij}$ are arbitrary constants.  \endproclaim

\demo{Proof} It remains only to explain why (3) is the general solution. Taking in (3) all
coefficients $C_{ij}$ but one equal to zero, we obtain $m_1+\dots+m_k=n$ partial solutions
$x^{m_i-l}e^{r_ix},\ l=1\,,\ldots,m_i,\ i=1\,,\ldots,k$. Their linear independence is well
known with various simple proofs by induction. \enddemo

{\bf Remark.} All numbers here may be both real and complex. In the last case one can get
real solutions, using the standard technique, such as Euler formula for exponent with complex
argument. \medskip

{\bf 3.} The same method can be used for solution of non-homogeneous equations as
well. Let $L[y]$ be the same operator as in Lemma 2 and let $y$ be a solution of an
equation $L[y] = f(x)$ with some integrable function $f$. Then the same proof as in
Lemma 2 gives that $(e^{-ax}y)^{(m)}=e^{-ax}f(x)$, so that
$$ y(x)= L^{-1}[f]=e^{ax}\int\cdots\int e^{-ax}f(x)\,dx=e^{ax}I_m[e^{-ax}f], \tag4 $$
where $I_m$ means an operator, inverse to $\frac{d^m}{dx^m}$. Of course, each
integral is defined here up to arbitrary constant, hence equality (4) presents
an explicit integral formula for the general solution of an equation $L[y] = f(x)$.

In the general case $L=L_1L_2\ldots L_k$ one has $L^{-1}=L_k^{-1}L_{k-1}^{-1}
\ldots L_1^{-1}$, so that the solution of the equation $L[y] = f(x)$ can be found by
iteration of formula (4):
$$ f_1=f,\quad f_{i+1}=e^{r_ix}I_{m_i}[e^{-r_ix}f_i],\quad i=1,\ldots,k,
   \qquad y=f_{k+1}(x). $$
If we ``drag" all appearing integration constants through all subsequent integrals,
we again obtain the general solution of the equation $L[y] = f(x)$. But in fact, all
these constants may be omitted, because we already have the general solution (3) for
the homogeneous equation and additionally need only one partial solution of the
equation $L[y] = f(x)$.

Let us consider, for instance, a case of $f(x)=e^{bx}Q_j(x)$, mostly studied in the higher
education. Here $Q_j(x)$ is a given polynomial of degree $j$ and $b$ may be either real or
complex, covering thus the trigonometric functions sine and cosine as well. As can be easily
verified via integrating by parts and omitting all integration constants, every integral
of $e^{cx}Q_j(x)$ with $c\neq0$ is a similar product $e^{cx}T_j(x)$ with the same $c$
and $j$, but with new coefficients of the polynomial $T_j$. Hence, if $b\neq r_1$, we get
$$ f_2(x)=e^{r_1x}\int\cdots\int e^{(b-r_1)x}Q_j(x)\,dx=e^{bx}T_j(x). \tag5 $$
Similar results will be obtained for all subsequent iterations until $b\neq r_i\,,\
i=2,\ldots,k$, and we get a final form of solution $y=e^{bx}S_j(x)$, where only
coefficients of the final polynomial $S_j(x)$ remain undetermined. Thus we can avoid
all long and tedious integrations, using only the method of undetermined coefficients.

Of course, the last result is well known as simply following from differential and integral
properties of the exponential and of the power functions. However, those arguments
cannot explain what to do when $b=r_i$, that is, when the exponent from the right-hand
side of a given equation coincides with one (at most!) root of the characteristic
polynomial. At the same time, the formula (5) immediately shows that in this case
the polynomial $Q_j(x)$ remains alone in all $m_i$ integrals, corresponding to the root
$r_i$, which increases its degree up to $j+m_i$. Note that this increasing concerns
all members of polynomial so that its smallest powers will disappear. At last, this
will look as replacing the final polynomial $S_j(x)$ by $x^{m_i}S_j(x)$.


\Refs
\ref\no1 \by W. E. Boyce and R. C. DiPrima \book Elementary Differential Equations 
and Boundary Value Problems \publ Willey and Sons, New York \yr 1997 \endref
\ref\no 2 \by E. Hille \book Lectures on Ordinary Differential Equations \publ 
Addison--Wesley, Reading \yr 1969 \endref

\endRefs
\end{document}